\documentclass[a4paper,twoside,11pt,reqno]{amsart}
\usepackage[T1]{fontenc}
\usepackage{courier}
\usepackage[scaled=0.92]{helvet}

\usepackage{amscd,amsmath,amssymb,amsthm,amsfonts,epsfig,graphics}

\usepackage[english]{babel}
\usepackage{amsthm,amsfonts,amssymb,amsmath}
\usepackage{hyperref}
\hypersetup{
colorlinks=true,
citecolor=red,
linkbordercolor={1 1 1},
citebordercolor={1 1 1},
}

\usepackage{mathrsfs}

\usepackage{xcolor} 
\theoremstyle{plain}
\newtheorem{satz}{Theorem}[section]
\newtheorem{prop}[satz]{Proposition}

\newtheorem{lem}[satz]{Lemma}

\theoremstyle{definition}

\newtheorem{rem}[satz]{Remark}

\newtheorem{hyp}[satz]{Hypothesis}

\newcommand{\mx}{\mbox}
\newcommand{\rw}{\rightarrow}
\newcommand{\de}{\displaystyle}
\newcommand{\ml}{\mathcal}

\newcommand{\pl}{\partial}

\newcommand{\x}{\times}

\newcommand{\beq}[1]{\begin{equation} \label{#1}}
\newcommand{\eeq}{\end{equation}}
\newcommand{\beqar}{\[ \begin{array}{rcl}}
\newcommand{\eeqar}{\end{array} \]}
\newcommand{\ue}{\underline}

\newcommand{\lie}[1]{ \mathcal{L}_{\chi^{(#1)}}}

\newcommand{\va}{\ue{\alpha}}

\newcommand{\jj}{^{(j)}}

\newcommand{\tilm}{\tilde{m}}

\providecommand{\ep}{\varepsilon}

\providecommand{\RR}{\mathbb{R}}
\providecommand{\CC}{\mathbb{C}}

\providecommand{\ZZ}{\mathbb{Z}}
\providecommand{\NN}{\mathbb{N}}
\newcommand{\snorm}[2]{\left| #1\right|_{#2}}
\newcommand{\norm}[2]{\left \lVert#1 \right\rVert_{#2}}

\newcommand{\Heq}[2]{\overset{\left(#1\right)}{\underset{}{#2}}}

\usepackage[utf8]{inputenc}
\hypersetup{pdfauthor={Fortunati, Wiggins},%
            pdftitle={Aperiodic Moser Theorem},%
            pdfsubject={Hamiltonian Mechanics},%
            pdfkeywords={Moser normal form, Aperiodic time dependence}%
}

\DeclareMathOperator{\id}{Id}

\makeatletter
\g@addto@macro{\endabstract}{\@setabstract}
\newcommand{\authorfootnotes}{\renewcommand\thefootnote{\@fnsymbol\c@footnote}}%
\makeatother

\begin{document}
\title[Normal forms \emph{\`{a} la Moser}...]{Normal forms \emph{\`{a} la Moser} for aperiodically time-dependent Hamiltonians in the vicinity of a hyperbolic equilibrium.}
\author{Alessandro Fortunati}
\thanks{This research was supported by ONR Grant No.~N00014-01-1-0769 and MINECO: ICMAT Severo Ochoa project SEV-2011-0087.}
\address{School of Mathematics, University of Bristol, Bristol BS8 1TW, United Kingdom}
\email{alessandro.fortunati@bristol.ac.uk}

\author{Stephen Wiggins}
\email{s.wiggins@bristol.ac.uk}

\maketitle

\begin{abstract}
The classical theorem of Moser, on the existence of a normal form in the neighbourhood of a hyperbolic equilibrium, is extended to a class of real-analytic Hamiltonians with aperiodically time-dependent perturbations. A stronger result is obtained in the case in which the perturbing function exhibits a time decay.
\end{abstract}

\section{Introduction}
The classical theorem of Moser, proven in \cite{moser56}, establishes the existence of a (convergent) normal form in a neighbourhood of a hyperbolic equilibrium of an area preserving map, either autonomous or periodically dependent on time. A result contained in \cite{cg94}, extends this result to the the flow of a priori unstable system in a neighbourhood of a partially hyperbolic torus, including in this way the quasiperiodic case. A concise description of the latter case can be found in \cite{g97}.\\
The aim of this paper is to show the existence of a normal form for Hamiltonians in the form (\ref{eq:ham}), i.e. real-analytic and non-autonomous perturbations of a hyperbolic equilibrium, for which the time dependence is not required to be periodic or quasiperiodic i.e. \emph{aperiodic}.\\
In the same spirit of the aperiodic version of the Kolmogorov theorem of \cite{fw14a}, which we use as a guideline (see also \cite{giortori}), the proof consists on the extension of the KAM approach of \cite{cg94} and \cite{g97}. 
Even in the original problem of Moser, despite the absence of ``genuine'' small divisors\footnote{This is a common feature with the ``non-purely hyperbolic'' case treated in \cite{giorlyap}.}, the well known property of \emph{superconvergence} of the KAM schemes, turns out to be of crucial importance in order to compensate the accumulation of ``artificial'' divisors generated by the Cauchy estimates. This feature is profitably used also in our case.\\
The treatment of the class of time-dependent homological equations, naturally arising in the normalization algorithm, has been improved with respect to \cite{fw14a}. Basically, the canonical transformation on which the single step of the mentioned algorithm is based, has the property to leave the time unchanged\footnote{This class of transformations was initially considered in \cite{gz92}.}. Hence, this can be interpreted as a family of canonical maps for which the time plays the role of ``parameter''. This allows to weaken the analyticity hypothesis for the time dependence leading to a remarkable simplification of the quantitative estimates.\\ 
The proof is carried out by using the formalism of the Lie series method developed by Giorgilli et al. (see e.g. \cite{gior02} and references therein).

\section{Preliminaries and statement of the result}
Let us consider the following Hamiltonian
\beq{eq:ham}
H(p,q,\eta,t)= \omega pq + \eta + F(p,q,t)
\mx{,}
\eeq
where $\omega \in (0,1]$, $(p,q,\eta) \in [-r,r]^2 \times \RR =: \ml{D}$ with $r>0$ and $t \in \RR^+:=[0,\infty)$. As usual, Hamiltonian (\ref{eq:ham}) is equivalent to the non-autonomous Hamiltonian $\ml{H}(p,q,t)=\omega p q + F(p,q,t)$ (which represents our original problem), by defining as $\eta$ the conjugate variable to $t$.\\
The function $F$ will be supposed to be real-analytic in $p$ and $q$ and such that, denoted as $f_{\ue{\alpha}}(t)$ its Taylor coefficients, one has $f_{\ue{\alpha}}(t)=0$ for all\footnote{It will be understood throughout the paper $\ue{\alpha} \in \NN^2$, denoting $|\va|:=\alpha_1+\alpha_2$.} $|\ue{\alpha}| \leq 2$, and all $t \in \RR^+$. Namely, the Taylor expansion of $F$ starts from the terms of degree $3$.\\  
The standard framework for the analysis, features the complexification of the domain $\ml{D}$ as follows.\\
Let $R \in (0,1/2]$ and define
$$
\ml{Q}_R:=\{(p,q) \in \CC^2 : |p|,|q| \leq  R \},\qquad \ml{S}_R:=\{\eta \in \CC : |\Im \eta| \leq  R \} \mx{,}
$$  
then set $\ml{D}_R:=\ml{Q}_R \times \ml{S}_{R} $. The perturbation $F$ will be supposed continuous on $\ml{Q}_{R}$ and holomorphic in the interior for all $t \in \RR^+$ (then $H$ is on $\ml{D}_{R} $) for some $R$. It will be sufficient to suppose that the real and imaginary parts of the complex valued functions $f_{\ue{\alpha}}(t)$ belong to $\ml{C}^1(\RR^+)$ for all $\ue{\alpha}$.\\
Given a function $G:\ml{Q}_R \times \RR^+ \rw \CC$, we consider the \emph{Taylor norm}  
\beq{eq:taylor}
\norm{G(p,q,t)}{R}:=\sum_{\ue{\alpha}} |g_{\ue{\alpha}}(t)|_+ R^{|\ue{\alpha}|} \mx{,} 
\eeq
where $|\cdot|_+:=\sup_{t \in \RR^+} |\cdot|$. Clearly $\snorm{G}{R}:=\sup_{\ml{Q}_R}|G|_+ \leq \norm{G}{R}$. We briefly recall the following standard result (which motivates the above described assumptions on $F$): if a function $G$ is continuous on $\ml{Q}_R$ and holomorphic in the interior, for all $t \in \RR^+$, one has $|g_{\ue{\alpha}}(t)|_+ \leq \snorm{G}{R}R^{-|\ue{\alpha}|}$. In particular, $\norm{G}{R'}<+\infty$ for all $R'<R$.\\
In the described setting the main result can be stated as follows
\begin{satz}[Aperiodic Moser '$56$]\label{thm:apermoser} Suppose that $1+\norm{F(p,q,t)}{R} =: M_F<\infty$. Then there exist $R_*,R_0$ with $0<R_*<R_0 \leq R^4$ and a family of canonical changes $\ml{M}:\ml{D}_{R_*} \rw \ml{D}_{R_0}$, analytic on $\ml{D}_{R_*}$ for all $t \in \RR^+$, casting the Hamiltonian (\ref{eq:ham}) in the \emph{time-dependent Moser normal form} 
\beq{eq:normalform}
H^{(\infty)}(p^{(\infty)},q^{(\infty)},\eta^{(\infty)},t)=J^{(\infty)}(x^{(\infty)},t)  + \eta^{(\infty)}  \mx{,}
\eeq
where $x:=pq$, $J^{(\infty)}(0,t)=0$ and $\pl_x J^{(\infty)}(0,t)=\omega$ for all $t \in \RR^+$.
\end{satz}
Exactly as in the classical Moser theorem, the quantity $x^{(\infty)}$ is a first integral, hence the flow associated to Hamiltonian (\ref{eq:normalform}) can be reduced to quadratures. In particular, one has
\[
p^{(\infty)}(t)=p^{(\infty)}(0) \exp(-\ml{A}(x^{(\infty)}(0),t)), \qquad 
q^{(\infty)}(t)=q^{(\infty)}(0) \exp(\ml{A}(x^{(\infty)}(0),t)), \qquad 
\]
where $\ml{A}(x,t):=\int_0^t \pl_x J^{(\infty)}(x,s)ds$.\\
The use of an additional ingredient leads to an even stronger result. Given $G:\ml{Q}_R \times \RR^+ \rw \CC$ we define as the ``time-dependent'' Taylor norm of $G$, the quantity  $\norm{G}{R;\RR^+}:=\sum_{\ue{\alpha}} |g_{\ue{\alpha}}(t)|R^{|\ue{\alpha}|}$, i.e. (\ref{eq:taylor}) in which $|\cdot|_+$ is replaced with $|\cdot|$. Now we introduce the next  
\begin{hyp}\label{hyp:decay}(Slow decay)
Suppose that there exist $M_F \in [1,+\infty)$ and $a>0$ such that 
\beq{eq:decay}
\norm{F(p,q,t)}{R;\RR^+} \leq M_F e^{-a t} \mx{,}
\eeq
for all $(p,q,t) \in \ml{Q}_R \times \RR^+$. 
\end{hyp}
In this way we are able to prove the following 
\begin{satz}[Strong Aperiodic Moser '$56$]\label{thm:strongmoser}
Under Hypothesis \ref{hyp:decay} it is possible to determine $0<\hat{R}_*<\hat{R}_0 \leq R^4$ and a family of canonical transformations $\ml{M}_{\ml{S}}$, analytic on $\ml{D}_{\hat{R}_*}$ for all $t \in \RR^+$, for which the Hamiltonian (\ref{eq:ham}) is transformed into the \emph{strong Moser normal form}
\beq{eq:strongnf}
\hat{H}^{(\infty)}(\hat{p}^{(\infty)},\hat{q}^{(\infty)},\hat{\eta}^{(\infty)},t)=\omega \hat{x}^{(\infty)} + \hat{\eta}^{(\infty)}  \mx{.}
\eeq
\end{satz}
The Hypothesis \ref{hyp:decay}, already used in \cite{fw14a}, turns out to be necessary in order to ensure the existence of certain improper integrals, which appear when dealing with time-dependent homological equations. As in the latter paper, this particular rate of decay is assumed only for simplicity of discussion. Similarly, we stress that no lower bounds are imposed on $a$ (except zero), in this way the time decay can be arbitrarily slow. The natural side-effect is that the estimates on the convergence radius of the normal form worsen as $a$ is smaller and smaller.\\
It should be stressed that, in both cases, the choice of $\omega$ in the interval $(0,1]$ is discussed as the ``interesting'' case. On the other hand, it is clear that the contribution of the time perturbation is smaller as $\omega$ increases\footnote{Namely, let $\mu:=O(\omega^{-1})$ and set $\hat{\omega}:=\mu \omega=O(1)$. Via a time rescaling $t=\mu \tau$,  problem (\ref{eq:ham}) is equivalent to the ``slowly time-dependent'' Hamiltonian $\hat{H}=\hat{\omega}pq+\mu \eta+ \mu F(q,u,\mu \tau)$.}. That is why, the case $\omega \geq 1$ can be treated with the same tools leading, in general, to easier estimates.\\
The proof of Theorem \ref{thm:apermoser} is (traditionally) achieved in two steps. In the first one (Sec. \ref{sec:formal}), a suitable normalization algorithm is constructed and discussed at a formal level. In the second part (Sec. \ref{sec:iterative}) the problem of its convergence is addressed, after having stated some tools of a technical nature (Sec. \ref{sec:technical}).\\
Proof of Theorem \ref{thm:strongmoser} is just a \emph{variazione sul tema}. The necessary modifications are outlined in Sec. \ref{sec:proofstrong}.  

%%%%%%%%%%%%%%%%%%%%%%%%%%%%%%%%%%%%%%%%%%%%%%%%%%%%%%%%%%%%%%%%%%%%%%%%%%%%%%%%%%
%% SECTION THE SECOND %%%%%%%%%%%%%%%%%%%%%%%%%%%%%%%%%%%%%%%%%%%%%%%%%%%%%%%%%%%%
%%%%%%%%%%%%%%%%%%%%%%%%%%%%%%%%%%%%%%%%%%%%%%%%%%%%%%%%%%%%%%%%%%%%%%%%%%%%%%%%%%

\section{The formal perturbative setting}\label{sec:formal}
The formal perturbative algorithm has the typical inductive structure. To start, we shall suppose that  Hamiltonian (\ref{eq:ham}) can be written at the $j-$th stage of the normalization process as
\beq{eq:hamricorsiva}
H^{(j)}(p,q,\eta,t)=\tilde{J}^{(j)}(x,t)+\eta+\tilde{F}^{(j)}(p,q,t) \mx{,}
\eeq
with $\tilde{F}^{(j)}$ at least of degree $3$ in $p,q$. It is immediate to realize that (\ref{eq:ham}) is in the form (\ref{eq:hamricorsiva}) so that we can set $H^{(0)}:=H$. Our aim is to construct a class of canonical transformations $\ml{M}_j$, parametrised by $t$, such that $H^{(j+1)}:=H^{(j)} \circ \ml{M}_j$ is still of the form (\ref{eq:hamricorsiva}). Roughly, the transformations $\ml{M}_j$ will be determined in such a way the ``mixed'' terms, i.e. of the form $p^{\alpha_1}q^{\alpha_2}$ with $\alpha_1 \neq \alpha_2$ contained in the perturbation, are ``gradually'' removed as $j$ increases, while the terms of the form $(pq)^n$ are progressively stored in $\tilde{J}^{(j)}$. This effect will be quantified in the next section, showing that the size of the  ``residual'' perturbation is asymptotic to zero, as $j \rw \infty$. Hence one sets 
\beq{eq:normaliz}
\ml{M}:=\lim_{j \rw \infty} \ml{M}_j \circ \ml{M}_{j-1} \circ \ldots \circ \ml{M}_0 \mx{,}
\eeq
so that, at least formally, $H^{(\infty)}=H \circ \ml{M}$.  \\
First of all we write 
\beq{eq:ftilde}
\tilde{F}^{(j)}(p,q,t)=\sum_{\substack{|\ue{\alpha}| \geq 3 \\
\alpha_1 \neq \alpha_2}
}\tilde{f}_{\ue{\alpha}}^{(j)}(t)p^{\alpha_1}q^{\alpha_2}+\sum_{k \geq 2} \tilde{f}_{\ue{k}}^{j}(pq)^k=:F^{(j)}(p,q,t)+\Delta^{(j)}(x,t) \mx{,}
\eeq
where $\ue{k}:=(k,k)$, then setting $J^{(j)}(x,t):=\tilde{J}^{(j)}(x,t)+\Delta^{(j)}(x,t)$, in such a way
\beq{eq:hamsis}
H^{(j)}(p,q,t)=J^{(j)}(x,t)+\eta+F^{(j)}(p,q,t) \mx{,}
\eeq
where $F^{(j)}$ contains only ``mixed'' terms.\\
Now we consider the action on $H^{(j)}$ of the transformation $\ml{M}_j$, which is defined by the the Lie series operator $\exp(\ml{L}_{\chi^{(j)}})=\id+\ml{L}_{\chi^{(j)}}+\
\sum_{s \geq 2}(1/s!) \ml{L}_{\chi^{(j)}}^s$. We recall that $\ml{L}_G F = \{F,G\}=F_q G_p+F_t G_{\eta} -F_p G_q - F_{\eta} G_t$, while $\chi^{(j)}=\chi\jj (p,q,t)$ is the (unknown) generating function.
Supposing that it is possible to determine it in such a way 
\beq{eq:homological}
\ml{L}_{\chi^{(j)}} (J^{(j)}(x,t)+\eta)+F^{(j)}(p,q,t)=0 \mx{,}
\eeq
one has that, by setting $\tilde{J}^{(j+1)}:=J^{(j)}$, and 
\beq{eq:fjputilde}
\tilde{F}^{(j+1)}:=\sum_{s \geq 1}\frac{1}{s!} \ml{L}_{\chi^{(j)}}^s F^{(j)} + \sum_{s \geq 2}\frac{1}{s!} \ml{L}_{\chi^{(j)}}^s (J^{(j)}+\eta) \mx{,}
\eeq
the transformed Hamiltonian $H^{(j+1)}:=\exp(\lie{j})H^{(j)}$ has exactly the form (\ref{eq:hamricorsiva}).\\Note that by (\ref{eq:homological}) and (\ref{eq:fjputilde}) 
\begin{align}
\tilde{F}^{(j+1)} & =  \de \sum_{s \geq 1} \frac{1}{s!}
\lie{j}^s \left[ F^{(j)}+\frac{1}{s+1} \lie{j} (J^{(j)}+\eta)\right] \nonumber \\
&= \de \sum_{s \geq 1} \frac{s}{(s+1)!} \lie{j}^s F^{(j)} \label{eq:truef} \mx{.}
\end{align}
Defining $g^{(j)}(x,t):=\pl_x J^{(j)}(x,t)$ one has that equation (\ref{eq:homological}) reads as 
\beq{eq:nwehom}
[g^{(j)}(x,t) \eth +\pl_t ]\chi^{(j)}(p,q,t)=F^{(j)} (p,q,t) \mx{,}
\eeq
having denoted $\eth:=q \pl_q-p \pl_p$. Taking into account of the expansion  $F^{(j)}=:\sum_{\substack{|\ue{\alpha}| \geq 3 \\ \alpha_1 \neq \alpha_2}} f_{\ue{\alpha}}^{(j)}(t) p^{\alpha_1} q^{\alpha_2}$, the solution of equation (\ref{eq:nwehom}) reads as 
\beq{eq:homscomposta}
\chi\jj(p,q,t)=\sum_{\ue{\alpha}} \ml{F}_{\ue{\alpha}}(x,t)p^{\alpha_1} q^{\alpha_2},\qquad
\ml{F}_{\ue{\alpha}}\jj (x,t):=e^{-\lambda A\jj (x,t)}\left[ \ml{F}_{\ue{\alpha},0}\jj (x) 
+\int_0^t e^{\lambda A\jj (x,s)} f_{\ue{\alpha}}\jj (s)ds
\right] \mx{.}
\eeq
where $A\jj (x,t):=\int_0^t g\jj (x,s)ds$, $\lambda:=\alpha_2-\alpha_1 \geq 1$ by hypothesis on $F^{(j)}$ and $\ml{F}_{\ue{\alpha},0}\jj (x)$ are functions to be determined. Clearly, we shall set $\ml{F}_{\ue{\alpha},0}\jj (x) \equiv 0$ for all $\ue{\alpha}$ such that $\alpha_1=\alpha_2$ and such that $f_{\ue{\alpha}}\jj (s) \equiv 0$ in such a way $\ml{F}_{\ue{\alpha}}\jj (x,t)$ are identically zero for those values.\\
It is evident that as $|\ue{\alpha}| \geq 3$ for by hypothesis on $F\jj$, the generating function $\chi\jj$ will be at least of degree $3$. This implies that,  by (\ref{eq:truef}), $\tilde{F}^{(j+1)}$ will be at least of degree $4$, in particular it will not contain terms of degree $2$. By hypothesis on $F \equiv F^{(0)}$ and by induction, this is true for all $j$, implying that 
$g^{(j)}(0,t)=\omega$ for all $t \geq 0$, i.e. $g^{(j)}$ has a strictly positive real part (by hypothesis on $\omega$), in a suitable neighbourhood of the origin and more precisely via a suitable choice of $R_0$. This will play a crucial role in our later arguments. The formal part is complete.
\begin{rem}
It is immediate to recognize the similarity between equation (\ref{eq:nwehom}) and those found in \cite{fw14a} and \cite{fw14c}. The main difference is the presence of the function $g^{(j)}(x,t)$ which requires a careful analysis about its variation on time, as anticipated above. 
\end{rem}

%%%%%%%%%%%%%%%%%%%%%%%%%%%%%%%%%%%%%%%%%%%%%%%%%%%%%%%%%%%%%%%%%%%%%%%%%%%%%%%%%%%%%%
%%%   SECTION THE THIRD   %%%%%%%%%%%%%%%%%%%%%%%%%%%%%%%%%%%%%%%%%%%%%%%%%%%%%%%%%%%%
%%%%%%%%%%%%%%%%%%%%%%%%%%%%%%%%%%%%%%%%%%%%%%%%%%%%%%%%%%%%%%%%%%%%%%%%%%%%%%%%%%%%%%

\section{Some preliminary results}\label{sec:technical}
\subsection{Bounds on the solutions of the homological equation}
First of all let us recall the following elementary equality, valid for all $\lambda \in [0,1)$, which will be repeatedly used in the follow
\beq{eq:useful}
\sum_{\va} \lambda^{|\va|}=\sum_{l \geq 0} (l+1)\lambda^l = (1-\lambda)^{-2} \mx{.}
\eeq
Then we state the next
\begin{prop}\label{prop:analitycity}
Suppose the existence of a positive constant $M\jj$ such that 
\beq{eq:smallf}
\norm{F\jj(p,q,t)}{R_j} \leq M\jj \mx{,}
\eeq
and that, for all $(x,t) \in \ml{Q}_{R_j} \times \RR^+$ one has
\begin{subequations}
\begin{align}
\Re g\jj (x,t) & \geq \omega/2 \label{eq:smallga} \mx{,}\\
|g\jj (x,t)| & \leq (3/2) \omega \label{eq:smallgb} \mx{.}
\end{align}
\end{subequations}
Then for all $\delta \in (0,1)$ the solution of (\ref{eq:nwehom}) satisfies
\beq{eq:chi}
\norm{\chi\jj (p,q,t)}{(1-\delta)R_j},\norm{\pl_t \chi \jj (p,q,t)}{(1-\delta)R_j} \leq \frac{4 M\jj}{\omega \delta^2} \mx{.}
\eeq
\end{prop}
\begin{rem}
Note that hypothesis (\ref{eq:smallga}) is essential as it is easy to find $g\jj (x,t)$ satisfying (\ref{eq:smallgb}) for which the solution of (\ref{eq:homscomposta}) is unbounded on $\RR^+$.
\end{rem}
The proof goes along the lines of a similar result contained in \cite{fw14a}, with the remarkable simplification due to the fact that now $t$ is purely real. The very minor drawback with respect to the ``analytic'' case treated in \cite{fw14a}, is that, in this case, the estimate of the time derivative does not follow directly from a Cauchy estimate.  
\proof
Recall that by hypothesis on $F$ one has $|f_{\ue{\alpha}}\jj (s)|_{+} \leq M\jj R_j^{-|\va|}$. If $\ue{\alpha}$ is such that $\lambda>0$, we shall set $\ml{F}_{\ue{\alpha},0}\jj (x)  \equiv 0$. By (\ref{eq:smallga}), we have that $\Re (A\jj (x,t)-A\jj (x,s))\geq \omega (t-s)/2$ on $\ml{Q}_{R_j}$, yielding
\beq{eq:lgz}
|\ml{F}_{\ue{\alpha}}\jj (x,t)| \leq M\jj  R_j^{-|\va|}e^{-\frac{\lambda \omega t}{2}} \int_0^t 
e^{-\frac{\lambda \omega s}{2}} ds  \leq \frac{2 M\jj }{\lambda \omega} R_j^{-|\va|} \mx{.}
\eeq
In the case $\lambda<0$, set $\lambda \rw - \lambda$ with $\lambda>0$, then we shall choose $\ml{F}_{\ue{\alpha},0}\jj(x):=-\int_{\RR^+}\exp(-\lambda A\jj (x,s)) f_{\va}(s)ds$. It is immediate to check that $|\ml{F}_{\ue{\alpha},0}\jj|<+\infty$ as, in particular, $\Re(A\jj (s))>\omega s /2$ by hypothesis. In such a way we get 
\beq{eq:lmz}
|\ml{F}_{\ue{\alpha}}\jj (x,t)| \leq M\jj R_j^{-|\va|}e^{\frac{\lambda \omega t}{2}} \int_t^{+\infty}  
e^{-\frac{\lambda \omega s}{2}} ds  \leq \frac{2 M\jj}{\lambda \omega} R_j^{-|\va|} \mx{.}
\eeq 
Hence,  by definition (\ref{eq:taylor}) and by  (\ref{eq:lgz}) and (\ref{eq:lmz}), for all $\lambda \in \ZZ \setminus \{0\}$
\beq{eq:normchi}
\norm{\chi\jj(p,q,t)}{(1-\delta)R_j} \leq \frac{2 M\jj }{|\lambda| \omega} \sum_{\va} (1-\delta)^{|\va|} \Heq{\ref{eq:useful}}{\leq} \frac{2 M\jj }{\omega \delta^2} \mx{,} 
\eeq
which implies the first part of (\ref{eq:chi}). \\
The second part of (\ref{eq:chi}) is straightforward from (\ref{eq:homscomposta}), bounds (\ref{eq:smallf}), (\ref{eq:lgz}), (\ref{eq:lmz}), and hypothesis (\ref{eq:smallgb}) then proceeding as in (\ref{eq:normchi}).
\endproof

\subsection{An estimate on the Lie operator}
This is a standard result in the works of A. Giorgilli et al., see e.g. \cite{giorlyap}. The statement recalled below, is adapted to the notational setting at hand
\begin{lem}\label{lem:lie}
Suppose that $\norm{\chi}{(1-\delta)R}$ and $\norm{G}{(1-\delta)R}$ are bounded for some $\delta \in (0,1/2)$. Then 
\beq{eq:lieest}
\norm{\ml{L}_{\chi}^s G}{(1-2 \delta)R} \leq s! (e^2 \delta^{-2} \norm{\chi}{(1-\delta)R})^s \norm{G}{(1-\delta)R}, \qquad \forall s \geq 1 \mx{.}
\eeq
\end{lem}
We shall also consider the case of bounded $\norm{G}{R}$, for which (\ref{eq:lieest}) clearly holds with the obvious replacement. It is evident that a sufficient condition for the convergence of the Lie operator $\exp(\ml{L}_{\chi})$ is that $e^2 \delta^{-2} \norm{\chi}{(1-\delta)R} \leq 1/2$.

%%%%%%%%%%%%%%%%%%%%%%%%%%%%%%%%%%%%%%%%%%%%%%%%%%%%%%%%%%%%%%%%%%%%%%%%%%%%%%%%%%%%%%
%%%   SECTION THE FOURTH  %%%%%%%%%%%%%%%%%%%%%%%%%%%%%%%%%%%%%%%%%%%%%%%%%%%%%%%%%%%%
%%%%%%%%%%%%%%%%%%%%%%%%%%%%%%%%%%%%%%%%%%%%%%%%%%%%%%%%%%%%%%%%%%%%%%%%%%%%%%%%%%%%%%

\section{Quantitative estimates}\label{sec:iterative}
\subsection{The iterative lemma}\label{subsec:iterative}
Let us consider a sequence $\{\ue{u}^{(j)}\}_{j \in \NN} \in [0,1]^5$ with $\ue{u}^{(0)}$ to be determined, where $\ue{u}^{(j)}:=(d_j,\ep_j,R_j,\tilde{m}_j,\tilde{M}_j)$. Let $\ue{u}_*:=(0,0,R_*,\tilde{m}_*,\tilde{M}_*)$ with $\omega/2  \leq \tilde{m}_* < \tilde{M}_* \leq (3/2) \omega$ and $R_*>0$ to be determined as well. Our aim is now to prove the next
\begin{lem} Suppose that for some $j \in \NN$, there exists $\ue{u}^{(j)}$ with $u_l^{(j)}>(u_*)_l$ for $l=1,\ldots,4$ and $\tilde{M}_j<\tilde{M}_*$, satisfying
\begin{align}
\norm{\tilde{F}^{(j)}(p,q,t)}{R_j} & \leq \ep_j \mx{,}
\label{eq:hypuno}\\
\Re \tilde{g}^{(j)}(x,t) & \geq \tilde{m}_j \mx{,}
\label{eq:hypdue} \\
|\tilde{g}\jj(x,t)| & \leq \tilde{M}_j \label{eq:hyptre} \mx{.}
\end{align}
for all $(x,t) \in \ml{Q}_{R_j} \times \RR^+$. 
Then, under the condition
\beq{eq:smallness}
\frac{4 e^2 \ep_j}{\omega R_*^2 d_j^6} \leq \frac{1}{2} \mx{,}
\eeq
it is possible to determine $u_l^{(j+1)} \in [(u_*)_l, u_l^{(j)})$ for $l=1,\ldots,4$ and $\tilde{M}_{j+1} \in (\tilde{M}_j,\tilde{M}_*]$ such that conditions (\ref{eq:hypuno}), (\ref{eq:hypdue}) and (\ref{eq:hyptre}) are satisfied by $\tilde{F}^{(j+1)}$ and $\tilde{g}^{(j+1)}$ as defined in Sec \ref{sec:formal}.
\end{lem}\label{lem:iterative}
The validity of (\ref{eq:hypdue}) and (\ref{eq:hyptre}) (compare with (\ref{eq:smallga}) and (\ref{eq:smallgb}))
with the above mentioned bounds on $\tilm_*$ and on $\tilde{M}_*$, is clearly related to the possibility of using Prop \ref{prop:analitycity} for all $j$ .
%It is clear the relation between (\ref{eq:hypdue}) and the possibility to use Prop. \ref{prop:time}. This yields a clear condition on $\ue{u}^*$ that we shall need to take into account, namely $\tilde{m}_j \geq  For this reason we shall set in the follows  hence we shall set a natural way to choose $\tilde{} it will be ensured that $\tilde{m}^*$
\proof
First of all, immediately from (\ref{eq:ftilde}) and (\ref{eq:hypuno}), it follows $\norm{F^{(j)}}{R_j}\leq \ep_j$. On the other hand, recall $g^{(j)}(x,t)=\tilde{g}\jj (x,t)+\pl_x \Delta\jj (x,t)$, where 
$\pl_x \Delta\jj (x,t) \equiv \sum_{k\geq 2} k f_{\ue{k}}\jj (t) x^{k-1} $, which implies
\[
\norm{\pl_x \Delta\jj (x,t)}{(1- 2 d_j)R_j} 
\leq \ep_j [(1- 2 d_j)R_j]^{-2} \sum_{k \geq 2} k(1- 2 d_j)^{2k} \leq  \ep_j (R_* d_j)^{-2} \mx{,} 
\]
hence on $\ml{Q}_{(1-2  d_j)R_j} \x \RR^+$
\beq{eq:gtildelow}
\Re g\jj (x,t) \Heq{\ref{eq:hypdue}}{\geq} \tilde{m}_j - \ep_j(R_* d_j)^{-2}=:m_j \mx{.}
\eeq
The last quantity is well defined as a consequence of the (stronger) condition (\ref{eq:smallness}), being  $\tilde{m}_j > \tilm_* \geq \omega/2$.  Similarly, $|g\jj (x,t)| \leq 
\tilde{M}_j + \ep_j (R_* d_j)^{-2}=:M_j
$.\\
From Lemma \ref{lem:lie} with $\delta=d_j$,  (\ref{eq:truef}), (\ref{eq:chi}) and (\ref{eq:hypuno}), under the \emph{convergence condition} guaranteed by (\ref{eq:smallness}) we get
\beq{eq:stimafpuno}
\norm{F^{(j+1)}(p,q,t)}{(1-2 d_j)R_j} \leq \ep_j \sum_{s \geq 1} \left(\frac{4e^2 \ep_j}{\omega d_j^4} \right)^s \leq 
\frac{8 e^2 \ep_j^2}{\omega d_j^4} \mx{.}
\eeq
Hence we shall set 
\beq{eq:recurrent}
\ep_{j+1}:=8 e^2 \omega^{-1} \ep_j^2 R_*^{-2} d_j^{-6}, \qquad R_{j+1}:=(1-2d_j)R_j, \qquad \tilde{m}_{j+1}:=m_j, \qquad \tilde{M}_{j+1}:=M_j \mx{,}
\eeq
in order to obtain the validity of (\ref{eq:hypuno}), (\ref{eq:hypdue}) and (\ref{eq:hyptre}) at the $j+1$-th step. The first of (\ref{eq:recurrent}) is the well known ``heart'' of the quadratic method.

\endproof
\subsection{Determination of the bounding sequences}\label{subsec:bounding}
Our aim is now to construct the sequence $\ue{u}\jj $ for all $j$ under the constraints (\ref{eq:recurrent}) and show that $\lim_{j \rw \infty} \ue{u}\jj  = \ue{u}_*$. The last step will be the determination of $\ue{u}_0$, completed in the next section. The procedure is analogous to \cite{fw14a}. We start by choosing, for all $j \geq 1$
\beq{eq:eps}
\ep_j:=\ep_0(j+1)^{-12} \mx{.}
\eeq
By substituting the latter into the first of (\ref{eq:recurrent}) we get $4 e^2 \ep_j/(\omega R_*^2 d_j^6) = 2^{-1}[(j+1)/(j+2)]^{12} \leq 1/2$, hence condition (\ref{eq:smallness}) holds for all $j \geq 0$. Similarly we get
\beq{eq:dj}
d_j=\left(\frac{8 e^2 \ep_0}{R_*^2 \omega}\right)^{\frac{1}{6}} \frac{(j+2)^2}{(j+1)^4} \mx{.}
\eeq
By supposing
\beq{eq:condepzeroone}
\ep_0 \leq 3^6 8^{-7} \pi^{-12} e^{-2} \omega R_*^2 \mx{,} 
\eeq
it is easy to see that
\beq{eq:summa}
\sum_{j \geq 0} d_j \leq 4 \left[ 8 e^2 \ep_0 (R_* \omega)^{-1} \right]^{\frac{1}{6}}  \sum_{j \geq 0} (j+1)^{-2} \leq 1/4 \mx{.}
\eeq
which implies, in particular, $d_j \leq 1/4$ for all $j \geq 0$ (essential for the correct definition of $R_{j+1}$). Condition (\ref{eq:condepzeroone}) will be obtained via a suitable choice of $R_0$ that will be addressed in Sec. \ref{subsec:trasf}.\\
By (\ref{eq:eps}), (\ref{eq:dj}), then by (\ref{eq:condepzeroone}), 
\[
(R_*)^{-2} \sum_{j \geq 0} \ep_j d_j^{-2} \leq
[8^{-1} \ep_0^2 \omega (R_*^2 e)^{-2} ]^{\frac{1}{3}} (\pi^2/6) < \omega/4 \mx{.}
\] 
Hence, comparing (\ref{eq:gtildelow}) with (\ref{eq:recurrent}), $\lim_{j \rw \infty} \tilde{m}_{j+1}=\tilde{m}_j-\ep_j(R_* d_j)^{-2} \geq \tilm_0 - (\omega/4)$. This implies that it is sufficient to set $\tilm_0=(3/4)\omega$ and $\tilm_* :=\omega/2$. Similarly we have $\lim_{j \rw \infty} \tilde{M}_j \leq \tilde{M}_*:=(3/2)\omega$ if $\tilde{M}_0:=(5/4)\omega$ is chosen.\\
As for $R_*$ we have $R_j:=R_0\prod_{l=0}^{j-1}(1-2 d_l)$.  
By writing $\log \prod_{l}(1-2 d_l)=\sum_{l} \log (1-2 d_l)$ and using (\ref{eq:dj}) under condition (\ref{eq:smallness}), we obtain\footnote{use inequality $\log(1-x) \geq -2 x \log 2$ (valid for all $x \in [0,1/2]$).} $\lim_{j \rw \infty} R_j \geq R_0/2=:R_*$. By replacing this value in (\ref{eq:dj}) and (\ref{eq:condepzeroone}), we see that $\ep_0$ and $d_0$ are determined once $R_0$ will be chosen.\\

\subsection{Transformation of variables and convergence of the scheme}\label{subsec:trasf}
For all $j \geq 0$ the transformation $\ml{M}_j:\ml{D}_{R_{j+1}} \rw \ml{D}_{R_j}$ acts on the variables as follows $(p\jj,q\jj,\eta\jj)=\lie{j}(p^{(j+1)},q^{(j+1)},\eta^{j+1})$, while $t$ is unchanged (as $\chi^{(j)}$ does not depend on $\eta$). Hence, by Lemma \ref{lem:lie}, then by the first of (\ref{eq:chi}) and condition (\ref{eq:smallness}), we get 
\beq{eq:boundp}
|p^{(j+1)}-p\jj| \leq \sum_{s \geq 1} (1/s!) \norm{\lie{j} p^{(j+1)} }{(1-2 d_j)R_j} \leq R_0^3 d_j^2/4 \mx{,} 
\eeq
analogously one obtains 
\beq{eq:boundq}
|q^{(j+1)}-q\jj| \leq R_0^3 d_j^2/4 \mx{.}
\eeq
As for $\eta$, write $\lie{j}^s \eta^{(j+1)} = -\lie{j}^{s-1} \pl_t \chi\jj$ then, similarly, by the second of (\ref{eq:chi}) 
\beq{eq:boundeta}
|\eta^{(j+1)}-\eta\jj| \leq 
\frac{4 \ep_j}{\omega d_j^2} \sum_{s \geq 1} \frac{(s-1)!}{s!} \left( \frac{4 e^2 \ep_j}{\omega d_j^4}\right)^{s-1} \leq \frac{d_j^2}{e^2}
\sum_{s \geq 1} \left( \frac{4 e^2 \ep_j}{\omega d_j^4}\right)^{s} \Heq{\ref{eq:smallness}}{\leq} \frac{R_0^2}{4 e^2} d_j^4 \mx{.}
\eeq
Our aim is now to determine the final value of $R_0$, by proceeding as follows. 
As $F$ is supposed to be analytic on $\ml{D}_{R}$, suppose $R_0 \leq R^4 \leq 1/16$. We have $|f_{\ue{\alpha}}|_+ \leq M_F R^{-|\ue{\alpha}|} \leq M_F R_0^{-|\ue{\alpha}|/4}$, hence (use (\ref{eq:useful}))
\[
\norm{F(p,q,t)}{R_0} \leq M_F \sum_{|\ue{\alpha}| \geq 3} R_0^{|\ue{\alpha}|/4} \leq 4 M_F R_0^{9/4}=:\ep_0 \mx{.} 
\] 
By substituting the latter into (\ref{eq:condepzeroone}) one gets
\beq{eq:rzeroone}
R_0 \leq \min\{(3^6 2^{-25} \pi^{-12} e^{-2} \omega M_F^{-1})^4,R^4\} \mx{.}
\eeq
Finally recall that $\tilde{J}^{(0)}=\omega x $ that is $\tilde{g}^{(0)}=\omega$. Hence, in order to guarantee that the choice of $\tilde{m}_0$ and of $\tilde{M}_0$ of Sec. \ref{subsec:bounding} is well defined, we need to show that $\norm{\pl_x \Delta^{(0)}}{R_0} \leq \omega/4$. Recall that $\Delta^{(0)}=\sum_{k \geq 2} f_{\ue{k}} x^k$, hence we get (use again the analyticity of $F$ on $\ml{D}_{R}$), $\norm{\pl_x \Delta^{(0)}}{R_0} \leq M_F \sum_{k \geq 2} k R_0^k \leq 8 M_F R_0^2$. It is immediate to realize that the latter is smaller than $\omega/4$, for all $\omega \in (0,1]$, under the condition (\ref{eq:rzeroone}). This completes the choice of $\ue{u}^{(0)}$.
\\
In conclusion, by using (\ref{eq:dj}) in   (\ref{eq:boundp}), (\ref{eq:boundq}) and (\ref{eq:boundeta}) we get, 
\[
\max\left\{\sum_{j \geq 0}|p^{(j+1)}-p\jj|,\sum_{j \geq 0}|q^{(j+1)}-q\jj|,\sum_{j \geq 0}|\eta^{(j+1)}-\eta\jj|\right\} \leq R_0/4 
\]
(we used $R_0 < 1 < e^2/\omega$, trivially from (\ref{eq:rzeroone})). Hence, by the Weierstra{\ss} theorem, the limit (\ref{eq:normaliz}) converges to a transformation, $\ml{M}:\ml{D}_{R_*} \rw \ml{D}_{R_0}$, which is analytic for all $t \in \RR^+$. Hence $(p^{(\infty)},q^{(\infty)},\eta^{(\infty)})$, denote the canonical variables on $\ml{D}_{R_*}$ (and $(p^{(0)},q^{(0)},\eta^{(0)}):=(p,q,\eta)$ those on $\ml{D}_{R_0}$) and the Hamiltonian $H^{(\infty)}$, formally defined after (\ref{eq:normaliz}), is an analytic function on $\ml{D}_{R_*}$ as well, and is in the desired Moser normal form.

\section{An outline of the proof of Theorem \ref{thm:strongmoser}}\label{sec:proofstrong}
In this section we describe the necessary modifications in the proof of Thm. \ref{thm:apermoser} in order to get its ``strong'' version. However, we stress that the crucial point is the following: if we suppose the existence of the integral $\int_{\RR^+} f_{\ue{\alpha}}\jj (t)dt$ (guaranteed by the exponential decay of $F\jj$), then (\ref{eq:homscomposta}) exists on $\RR^+$ also for $\lambda=0$ i.e. the r.h.s. of the homological equation can contain also terms with $\alpha_1=\alpha_2$. 
\subsubsection*{Formal scheme}
The definition of $\tilde{J}\jj$ and of $\tilde{F}\jj$ is not necessary, we suppose that $H\jj$ is directly of the form 
\beq{eq:hamricstrong}
H\jj=\omega p q + \eta + F\jj (p,q,t) \mx{.}
\eeq
The initial Hamiltonian is exactly of the form above, so we can set $H^{(0)}:=H$. Suppose that $\chi\jj$ is chosen in a way to satisfy the homological equation
\beq{eq:stronghom}
\lie{j}(\omega p q + \eta )+F\jj=0 \mx{,}
\eeq
it is sufficient to define
\beq{eq:fjponestrong}
F^{(j+1)}:=\sum_{s \geq 1} \frac{s}{(s+1)!}\lie{j} F\jj \mx{,}
\eeq
in order to have $H^{(j+1)}$ of the form (\ref{eq:hamricstrong}). By expanding $\chi\jj=\sum_{\va} c_{\va}(t)p^{\alpha_1} q^{\alpha_2}$ and $F\jj$ as well\footnote{Note that in this case the Taylor expansion of $F\jj$ will contain also terms with $\alpha_1=\alpha_2$.}, we get this time, for all $\ue{\alpha}$
\beq{eq:homologicalstrong}
\dot{c}_{\ue{\alpha}}\jj (t)+\hat{\lambda} c_{\ue{\alpha}}\jj (t)=f_{\ue{\alpha}}\jj (t)\mx{,}
\eeq 
with $\hat{\lambda}:=\omega (\alpha_1-\alpha_2)$ \emph{purely real}.
\subsubsection*{Bounds on the homological equation}
The easy structure of eq. (\ref{eq:homologicalstrong}) simplifies remarkably the proof of the equivalent of Prop. \ref{prop:analitycity}, which states, in this case, as follows
\begin{prop} Suppose that there exists $M\jj >0$ such that $\norm{F\jj(p,q,t)}{R_j,\RR^+} \leq M\jj \exp(-a t)$. Then for all $\delta \in (0,1)$ the solution of (\ref{eq:homologicalstrong}) satisfies
\beq{eq:estimatestrong}
\norm{\chi\jj(p,q,t)}{(1-\delta)R_j;\RR^+},\norm{\pl_t \chi\jj(p,q,t)}{(1-\delta)R_j;\RR^+} \leq 4 M\jj a^{-1} \delta^{-3} \mx{.}
\eeq 
\end{prop}
\proof(Sketch)
If $\ue{\alpha}$ is such that $\hat{\lambda}>0$ then choose $c_{\ue{\alpha}}\jj (0)=0$. In this way\footnote{We are using $\exp(\hat{\lambda}(s-t)) \leq 1$ (and a similar bound in the case $\hat{\lambda} \leq 0$). These bounds, used for simplicity and sufficient for our purposes, hide the property $|c_{\ue{\alpha}}\jj(t)| \rw 0$ as $t \rw +\infty$ i.e. the (reasonable) phenomenon for which $\chi\jj$ is asymptotically vanishing for all $j$, that is, the canonical transformation $\ml{M}_{\ml{S}}$ reduces to the identity.} 
\beq{eq:estimatecstrong}
|c_{\ue{\alpha}}\jj(t)| \leq M\jj R_j^{-|\ue{\alpha}|}
\int_0^t e^{\hat{\lambda}(s-t)} e^{-a s} ds \leq M\jj a^{-1} R_j^{-|\ue{\alpha}|}  \mx{.}
\eeq
If $\hat{\lambda} \leq 0$ set $\hat{\lambda} \rw -\hat{\lambda}$ with $\hat{\lambda} \geq 0$ and choose $c_{\ue{\alpha}}\jj (0):=-\int_0^{+\infty} \exp(-\hat{\lambda}s) f_{\ue{\alpha}}\jj (s)ds$. A similar procedure yields the same estimate as (\ref{eq:estimatecstrong}) and then $\norm{\chi\jj(p,q,t)}{(1-\delta)R_j;\RR^+} \leq M a^{-1} \delta^{-2}$. By using the obtained estimates and (\ref{eq:homologicalstrong}), one gets the second of (\ref{eq:estimatestrong}).  
\endproof
\subsubsection*{Quantitative part}
We define now $\ue{\hat{u}}\jj :=(d_j,\ep_j,R_j)$, with $\ue{\hat{u}}_*:=(0,0,\hat{R}_*)$. Statement of Lemma \ref{lem:iterative} modifies as follows
\begin{lem}
Suppose that for some $j \in \NN$, there exists $\ue{\hat{u}}^{(j)}$ with $\hat{u}_l^{(j)}>(\hat{u}_*)_l$ for $l=1,2,3$, satisfying
\beq{eq:hypunobis}
\norm{\tilde{F}^{(j)}(p,q,t)}{R_j;\RR^+}  \leq \ep_j e^{-a t}\mx{,}
\eeq 
for all $(p,q,t) \in \ml{Q}_{R_j} \times \RR^+$. 
Then, under the condition
\beq{eq:smallnessstrong}
\frac{4 e^2 \ep_j}{\omega a \hat{R}_*^2 d_j^6} \leq \frac{1}{2} \mx{,}
\eeq
it is possible to determine $\hat{u}_l^{(j+1)} \in [(\hat{u}_*)_l, \hat{u}_l^{(j)})$ for $l=1,2,3$ such that (\ref{eq:hypunobis}) is satisfied by $F^{(j+1)}$ as defined in (\ref{eq:fjponestrong}).
\end{lem}
The proof of this Lemma and of the rest of the Theorem is straightforward \emph{mutatis mutandis}. We only mention that condition (\ref{eq:condepzeroone}) is replaced by $\ep_0 \leq 3^6 8^{-7} e^{-2} \pi^{-12} a \omega \hat{R}_*^2$, implying 
\[\hat{R}_0 \leq \min \{(3^6 2^{-25} \pi^{-12} e^{-2} \omega a M_F^{-1})^4,\hat{R}^4\}
\]
i.e. $\hat{R}_0 \sim a^4$ as $a \rw 0$ as announced after the statement of Thm. \ref{thm:strongmoser}. 
\subsection*{Acknowledgements} We are indebted with Prof. Antonio Giorgilli for some enlightening comments. The first author is also grateful to Prof. Luigi Chierchia for useful discussions on his cited work.

\bibliographystyle{alpha}
\bibliography{MoserNF.bib}

\newcommand{\sort}[1]{}
\begin{thebibliography}{FW14b}

\bibitem[CG94]{cg94}
L.~Chierchia and G.~Gallavotti.
\newblock Drift and diffusion in phase space.
\newblock {\em Ann. Inst. H. Poincar\'e Phys. Th\'eor.}, 60(1):144, 1994.

\bibitem[FW14a]{fw14a}
A.~Fortunati and S.~Wiggins.
\newblock Persistence of {D}iophantine flows for quadratic nearly integrable
  {H}amiltonians under slowly decaying aperiodic time dependence.
\newblock {\em Regul. Chaotic Dyn.}, 19(5):586--600, {\sort{a}} 2014.

\bibitem[FW14b]{fw14c}
A.~Fortunati and S.~Wiggins.
\newblock {A Kolmogorov theorem for nearly-integrable Poisson systems with
  asymptotically decaying time-dependent perturbation.}
\newblock \url{http://arxiv.org/abs/1409.0430}, {\sort{b}} 2014.

\bibitem[Gal97]{g97}
G.~Gallavotti.
\newblock Hamilton-{J}acobi's equation and {A}rnold's diffusion near invariant
  tori in a priori unstable isochronous systems.
\newblock {\em Rend. Sem. Mat. Univ. Politec. Torino}, 55(4):291--302 (1999),
  1997.
\newblock Jacobian conjecture and dynamical systems (Torino, 1997).

\bibitem[Gioa]{giorlyap}
A.~Giorgilli.
\newblock On a {T}heorem of {L}yapounov.
\newblock {\em Rendiconti dell'Istituto Lombardo Accademia di Scienze e
  Lettere, Classe di Scienze Matematiche e Naturali}, In Press.

\bibitem[Giob]{giortori}
A.~Giorgilli.
\newblock Persistence of invariant tori.
\newblock \url{http://www.mat.unimi.it/users/antonio/hamsys/hamsys.html}.

\bibitem[Gio03]{gior02}
A.~Giorgilli.
\newblock Exponential stability of {H}amiltonian systems.
\newblock In {\em Dynamical systems. {P}art {I}}, Pubbl. Cent. Ric. Mat. Ennio
  Giorgi, pages 87--198. Scuola Norm. Sup., Pisa, 2003.

\bibitem[GZ92]{gz92}
A.~Giorgilli and E.~Zehnder.
\newblock Exponential stability for time dependent potentials.
\newblock {\em Z. Angew. Math. Phys.}, 43(5):827--855, 1992.

\bibitem[Mos56]{moser56}
J.~Moser.
\newblock The analytic invariants of an area-preserving mapping near a
  hyperbolic fixed point.
\newblock {\em Comm. Pure Appl. Math.}, 9:673--692, 1956.

\end{thebibliography}

\end{document}